\documentclass{amsart}
\usepackage{bbm}
\usepackage{graphicx}
\usepackage{verbatim}
\usepackage{stackrel}
\usepackage{color}
\graphicspath{{Fig/}}
\usepackage{amsmath}
\usepackage{amssymb}
\usepackage{array}
\usepackage{mathtools}
\usepackage{xcolor}
\usepackage{mathdots}
\usepackage{pgf}
\usepackage{tikz}
\usetikzlibrary{shapes.arrows, fadings}
\usetikzlibrary{arrows}
\usetikzlibrary{matrix}
\usetikzlibrary{shapes,intersections}

\usepackage[latin1]{inputenc}

\newtheoremstyle%
 {bluethm}%
 {}{}%
 {\color{blue}\itshape}
 {}%
 {\color{blue}\bfseries}%
 {\color{blue}.}%
 { }{}

 \newtheoremstyle%
 {redthm}%
 {}{}%
 {\color{red}\itshape}
 {}%
 {\color{red}\bfseries}%
 {\color{red}.}%
 { }{}

\newtheorem{theorem}{Theorem}
\newtheorem{prop}{Proposition}
\newtheorem{lemma}{Lemma}
\newtheorem{defi}{Definition}
\newtheorem{assumption}{Assumption}[part]
\theoremstyle{definition}

\newtheorem{remark}{Remark}
\newtheorem{example}{Example}

\def\N{{\mathbb N}}
\def\R{{\mathbb R}}

\def\P{{\mathbb P}}
\def\E{{\mathbb E}}

\newcommand{\diff}{\mathop{}\mathopen{}\mathrm{d}}

\newcommand\ind[1]{\mathbbm{1}_{\left\{#1\right\}}}

\def\cal{\mathcal}
\def\eps{\varepsilon}

\title[ Condensation in reversible Coag-Frag models]{On the Condensation and fluctuations in reversible coagulation-fragmentation models}

\author{Wen Sun}
\address        {School of Mathematical Sciences, University of Science and Technology of China, Jinzhai 96, 230026 Hefei}

\date{\today}

\begin{document}

\maketitle
We study the condensation phenomenon for the invariant measures of the mean-field model of reversible coagulation-fragmentation processes conditioned to a supercritical density of particles. It is  shown that when the parameters of the associated balance equation satisfy a subexponential tail condition, there is one single giant particle that corresponds to the missing mass in the macroscopic limit. We also show that in this case, the rest of the particles are asymptotically \emph{i.i.d.} according  to the normalized equilibrium state of the limit hydrodynamic differential equation. Conditions for the normal fluctuations and the $\alpha$-stable fluctuations around the condensed mass are given. We obtain the large deviation principle  for the empirical measure of the masses of the particles at equilibrium as well.

\bigskip\bigskip

\hrule

\vspace{-0.4cm}

\tableofcontents

\vspace{-0.9cm}

\hrule

\bigskip

\section{Introduction and main results}

We study the condensation phenomenon in the mean-field models of random reversible coagulation-fragmentation particle systems.

\subsection{The coagulation-fragmentation particle systems}
We recall that  coagulation-fragmentation particle systems concern the dynamics of the masses of all particles due to two types of interactions:
\begin{enumerate}
\item any two particles  can coagulate into a new particle;
\item any particle can break into two smaller particles,
\end{enumerate}
following the rule of mass conservation. Given two symmetric functions $a,b:(\N^*)^2=\{1,2,\dots\}^2\mapsto \R^+$, the dynamics of the system is described as follows. If there are two particle of mass $i$ and $j$ in the system, then independently of all the other particles, this two particles are coagulated after an exponential waiting time with rate $a(i,j)$; for all particle of mass $i+j$ in this system, independently of all the others, it can break into two particles of masses $i$ and $j$ after an exponential waiting time with rate $b(i,j)$. For a given volume $V\in \R^+$ and total mass $M\in\N^*$, a mean-field coagulation-fragmentation process can be rigorously defined as a  Markov process on the state space
\[
\cal{S}^{M}:=\left\{\eta\in \N^M\bigg|\sum_{i=1}^M i\eta_i=M\right\},
\]
with an infinitesimal generator
\begin{multline}\label{cfgen}
  \cal{L}^{V,M}f(\eta)=\sum_{1\le i\le j\le M}\frac{a(i,j)}{V}\eta_i\eta_j\left(f(\eta^{i,j}_+)-f(\eta)\right)\\
  +\sum_{1\le i\le j\le M}{b(i,j)}\eta_{i+j}\left(f(\eta^{i,j}_-)-f(\eta)\right),
\end{multline}
where
\[
\eta^{i,j}_+=\eta+\delta_{i+j}-\delta_i-\delta_j,
\]
and
\[
\eta^{i,j}_+=\eta-\delta_{i+j}+\delta_i+\delta_j.
\]

If there exists a sequence of positive numbers $(Q_r,r\ge 1)$
such that for all $i,j\in\N^*$, the following details balanced equation
\begin{equation}\label{balance}
a(i,j)Q_iQ_j=Q_{i+j}b(i,j)
\end{equation}
holds, then we call this particle system with generator $\cal{L}^{V,M}$ a reversible coagulation-fragmentation system. Under this assumption, the unique invariant measure of the system is given by
\begin{equation}\label{inme}
\pi^{V,M}(\eta)=\frac{1}{Z^{V,M}}\prod_{i=1}^M\frac{\left(V Q_i\right)^{\eta_i}}{\eta_i!},
\end{equation}
where $Z^{V,M}$ is the normalisation constant. For more about the reversible particle system and its invariant measure, we refer to the book Kelly~\cite{kelly}. We remark here that if a sequence $(Q_r,r\ge 1)$ satisfies the balance equation~\eqref{balance}, then for any $\phi>0$, the sequence $(Q_r\phi^r,r\ge 1)$ also satisfies the balance equation.

\begin{example}[The Becker-D\"oring coagulation-fragmentation model] If one particle can only coagulates with a monomer, which is a particle with mass $1$, and one particle can only breaks into one smaller particle and one monomer, then this system is called a Becker-D\"oring model. We claim that when the reaction rates are positive, then the system is reversible. To see this, we let
  \[
a_r:=a(r,1),\qquad b_{r+1}:=b(r,1),
\]
and define a sequence $(Q_r)$ by
\[
Q_1=1;\qquad \frac{Q_{r+1}}{Q_r}=\frac{a_r}{b_{r+1}},\qquad \forall r\ge1,
\]
then $(Q_r)$ satisfies the detailed balance equation~\eqref{balance}. For more about the Becker-D\"oring model, we refer to the survey Hingant et al.~\cite{hingant2017deterministic}.
\end{example}
\subsection{The condensation in the hydrodynamic coagulation-fragmentation equation}
Let $\eta^V(t)=(\eta^V_1(t),\eta^V_2(t)\dots,\eta^V_M(t))$ be the Markov jump process with the infinitesimal generator~\eqref{cfgen} described above. It has been proved by Jeon~\cite{jeon1998existence}, that when the reaction rates $a(i,j)$ and $b(i,j)$ are with order $o(ij)$ and the density has a finite limit  $\lim_{V\to\infty}M/V=\rho<\infty$, then on any finite time interval $[0,T]$, the process $\eta^V(t)/V$ converges weakly in probability to the  coagulation-fragmentation Equation,
\[
(\eta^V(t)/V,t\in[0,T])\Rightarrow (c(t),t\in[0,T]),
\]
where for all $r\ge 1$,
 \begin{multline}\label{cfeq}
    {c}_r'(t)=\frac{1}{2}\sum_{i=1}^{r-1}a(i,r-i)c_i(t)c_{r-i}(t)-c_r(t)\sum_{i=1}^\infty a(r,i)c_i(t)\\
    +\sum_{i=1}^\infty b(i,r)c_{i+r}(t)-\frac{1}{2}\sum_{i=1}^{r-1}b(r-i,i)c_r(t).
 \end{multline}
In the paper~\cite{ldpsun}, we have proved the associated pathwise large deviation principle under a more general setting. For the Becker-D\"oring model, a functional central limit theorem has been established by the author in~\cite{sun2018functional}.

The Becker-D\"oring coagulation-fragmentation equation has been  studied  by Ball et al.~\cite{ball}. It has been show that, if  there exists a finite $\varphi(\rho)$ such that
\[
\sum_{r\ge 1}rQ_r\varphi(\rho)^r=\rho
\]
then this equation has exactly one equilibrium state $c^\rho$ with density $\rho$, 
\[
c^\rho_r=Q_r\varphi(\rho)^r,\qquad \forall r\ge 1.
\]
Moreover, under additional conditions, the solution of the  coagulation-fragmentation Equation~\eqref{cfeq} $c(t)$ converges in the pointwise topology to this equilibrium solution $c^\rho$ as $t\to\infty$. When such a  finite $\varphi(\rho)$ does not exists, let $\rho_c$ be the largest value such that the function $\sum_{r\ge 1}rQ_r\varphi(x)^r=x$ has a root. It is clear $\rho_c<\rho$. Then the result in Ball et al.~\cite{ball} tells us that,  under additional conditions, the solution  $c(t)$ converges in the pointwise topology to this critical equilibrium solution $c^{\rho_c}$ as $t\to\infty$. That is to say, in the hydrodynamic regime, if $\rho>\rho_c$, then there is some mass missing as time $t\to\infty$, which corresponds to the condensation phenomenon in the Becker-D\"oring models. We remark that the results in  Ball et al.~\cite{ball} can be extended to the general reversible coagulation-fragmentation equation without much efforts.

\subsection{Main contribution} In this paper, we study the condensation phenomenon in the stochastic coagulation-fragmentation model. We show that under the invariant measure $\pi^{V,M}$, when the density $M/V$ converges to a finite $\rho$, then the scaled random sequence $(\eta_r/V,r\ge 1)$ converges in $\{\ell\in\R_+^\infty|\sum_r r\ell_r\le \rho\}$ with the pointwise topology to the equilibrium solution $c^{\rho\land \rho_c}$ by establishing a large deviation principle with a good rate function
\[
\sum_{r=1}^\infty\left(\ell_r\log\left(\frac{\ell_r}{c_r^{\rho\land \rho_c}}\right)-\ell_r+c_r^{\rho\land \rho_c}\right)-\left(\rho-\sum_{r=1}^\infty r\ell_r\right)\log \frac{\varphi(\rho)}{\phi_c},
\]
where $\phi_c=\lim_{r\to\infty}Q_r^{1/r}<\infty$ and $\varphi(\rho)=\phi_c$ if $\rho\ge \rho_c$. We can see that the first term is the relative entropy between the measure $(\ell_r,r\ge 1)$ and  $(c_r^{\rho\land \rho_c},r\ge 1)$ on positive integers. The second term describes the difference of total mass. Our result build a relation between the condensation phenomenons in the microscopic (stochastic) and the macroscopic (hydrodynamic) reversible coagulation-fragmentation models.

Later, we give an sufficient condition on the sequence $(Q_r,r\ge 1)$ and prove that, under this condition, there is only one giant particle appears under the invariant measure $\pi^{V,M}$ when $M$ sufficient large  in the limit. The rest of the particles are asymptotically \emph{i.i.d.} distributed. Then under additional assumptions, we can establish a central limit theorem for this giant particle. For example, when $\sum_r r^2Q_r\phi_c^r<\infty$, then the fluctuation is approximately a Gaussian with variance $\sum_r r^2Q_r$; when $\sum_r Q_r\phi_c^r\approx 1/r^{\alpha}$ for some $\alpha>0$, then the fluctuation is approximately a $\alpha$-stable  random variable. We should emphasise that this assumption can go beyond the hydrodynamic regime -- the limit density $M/V$ can be infinite large. When the density $M/V$ converges  to a finite $ \rho>\rho_c$, then the mass of the giant particle is approximately $V(\rho-\rho_c)$ which corresponds to the missing mass in the equilibrium state of the coagulation-fragmentation equation~\eqref{cfeq}. The fluctuation around $\eta_r$ are also given under the invariant measure $\pi^{V,M}$. Moreover, it is  shown that for any $r\neq k$, $\eta_r$ and $\eta_k$ are asymptotically independent. It is different from the finite marginal central limit theorems in non mean-field coagulation-fragmentation model proved by Durrett et al.~\cite{Durrett99}, where the covariance between $\eta_r$ and $\eta_k$ are given.

\subsection{Relations with other models}

\subsubsection{The appearance of a single giant `particle' in random walks, Zero-Range models and Erd\"os-R\'enyi graphs} It turns out that the condensation phenomenon in our stochastic coagulation-fragmentation model is related to many famous phase transition in the mass conserving models.

For example, in a random walk $(X_1,X_2,\dots, X_n)$ with a constraint $\{\sum_{i=1}^nX_i=m\}$, if the distribution of $X$ belongs to a subclass of the sub-exponential distribution, then Denisov et al.~\cite{denisov} prove that when $m$ is sufficient large, there will be only one giant jump appears in this random walk. There are huge literature in this domain. We refer to  Nagaev~\cite{MR0282396}, Tka\v{c}uk~\cite{MR0368115}, Baltrunas~\cite{MR1429811}, Kl\"{u}ppelber et al.~\cite{kumi} and also the references in  Denisov et al.~\cite{denisov}.
The Zero-Range model also owns its condensation phase transition. When the jump rate function satisfies a Power law condition, if the total mass of the model is sufficient large, then a giant site emerges (Gro{\ss}kinsky et al.~\cite{gross}). The relation between the condensation in Zero-Range models and the giant jump in subexponential random walk has been shown by Armend\'ariz et al.~\cite{armen, armendariz}.
See also Jeon et al.~\cite{jeon2}.

It will be shown later, the condensation in the stochastic coagulation-fragmentation model can be represented by the appearance of the giant jump in a compound Poisson process under certain subexponential condition. In the paper by the author~\cite{sunmdp}, we apply a similar method to study the sparse Erd\"os-R\'enyi random graph. We are able to prove the moderate deviation principle for the empirical measures of the sizes of the connected components and also for the giant cluster. Despite that the component Poisson process associated with the random graph is more complicated, the reason for the emergence of a giant `particle' in the sparse graph and the in the coagulation-fragmentation model are similar. We remark that this similarity is not that surprising since the distribution of the sparse random graph is exactly the distribution of the Smoluchowski's pure coagulation model with a product kernel. See Aldous' survey~\cite{aldous1999deterministic} for example.

\subsubsection{Cluster expansion} The statistical mechanics of particle systems can be studied by Mayer's cluster expansions~\cite{Mayer}. The partition function can be given by
\begin{equation}\label{part}
\sum_{i\eta_i=M}\prod_{i=1}^M\frac{1}{\eta_i!}\left(\frac{V b_i(V,T)}{\lambda^d}\right)^{\eta_i},
\end{equation}
where $\lambda$ is the thermal wavelength and $b_k(V,T)$ are the cluster integrals depends on the volume $V$ and temperature $T$. In many cases, the cluster integrals are positive. For example, in an ideal Bose gas model, we have that for any finite $k$, and volume $V$ sufficient large,
\[
\frac{b_k(V,T)}{\lambda^d}\approx \left(\frac{1}{4\pi\beta}\right)^{d/2}\frac{1}{k^{1+d/2}},\]
where $\beta=1/T$.
We notice that the partition function~\eqref{part} has a similar shape as the normalizer of our invariant measure $\pi^{V,M}$. Indeed, the papers by Adams~\cite{adams2007large,adams2008large} calculate the free energy and prove a large deviation principle for $(\eta_1/V,\eta_2/V,\dots)$ under $\pi^{V,M}$ with
\[
Q_k=\left(\frac{1}{4\pi\beta}\right)^{d/2}\frac{1}{k^{1+d/2}},\]
to study the Bose-Einstein condensation in the ideal Bose gas. Similar large deviation principles have been established in various models, see for example Andreis et al.~\cite{andreis2021large} and also Jansen et al.~\cite{jansen2015large}. While their proofs rely on the most probable state under a certain measure, our method relies on the behaviours of a conditional compound Poisson process and gives a less technical proof.

\subsubsection{Random permutation with large cycle weights}
From the paper Betz et al.~\cite{betz}, we find that $\pi^{V,M}$ is also the probability of the occupation numbers in the random permutation model of $M$ elements and large cycle weights $(rQ_rV,r\ge 1)$. Our results show the emergence of a giant cycle in this random permutation. In the some partition models for the integer $n$, it has been shown that large cycles appears with an order $\sqrt{n}$. Moreover, the shape of large cycles are captured by Vershik~\cite{vershik} and see also Dembo et al.~\cite{dvz}. It would be interesting for us to look at the relation between the coagulation-fragmentation model and the random partition model in the future.

\subsection*{Outline of the paper}
This paper is organised as follows. We firstly state the notations and our main results in Section~\ref{re}. In Section~\ref{sec:ldp} we prove the large deviation principle (Theorem~\ref{main:ldp}) for the scaled random variable $(\eta_1/V,\eta_2/V,\dots)$ under the invariant measure $\pi^{V,M}$.
In Section~\ref{sec:cond} we introduce a compound Poisson process with a random summation condition that gives a new representation of the measure $\pi^{V,M}$. We also prove the appearance of the giant jump and the equivalence of ensemble properties (Proposition~\ref{poisum}) in this conditional compound Poisson process. With the help of it, we are able to prove the  condensation phenomenon (Theorem~\ref{main:cond}) and the associated central limit theorems (Theorem~\ref{main:clt}) in Section~\ref{sec:cpp}.

\subsection*{Notations}
We use the notation $\N=\{0,1,2,\dots\}$, $\N^*=\{1,2,3,\dots\}$, $\R_+=[0,\infty)$, $\overline{\R}_+=[0,\infty]$. We say a function $L$ is slowly varying, if for all fixed $t>0$,
\[
\lim_{x\to\infty}\frac{L(tx)}{L(x)}=1.
\]

\section{Basic assumptions and main results}\label{re}

In this section we  introduce the notations and assumptions used throughout this paper. We  state our main result on the large deviation principle under the invariant measure $\pi^{V,M}$ when the density $M/V\to\rho$. We give a sufficient condition on $(Q_r)$, under which, when the density $M/V$ is sufficient large, only one giant particle emerges that represents the condensation. With additional assumptions, the fluctuation around the size of this giant particle can be depicted.
\begin{defi}
For any $(Q_r)$, we define a mapping $F:\R_+\mapsto \overline{\R}_+$ by
  \[
F(\phi):=\sum_{r=1}^\infty rQ_r\phi^r.
\]
Let $\phi_c$ be the radius of convergence.
Define and suppose
\[
\rho_c:=\sup_{0\le \phi<\phi_c}F(\phi)\in\overline{\R}_+.
\]
Clearly, $F$ is smooth and strictly increasing on $(0,\phi_c)$ and $\lim_{\phi\to\phi_c-}F(\phi)=\rho_c$.
For $\rho<\infty$ and $0\le \rho\le \rho_c$, let $\varphi(\rho)$ be the unique root of $F(\phi)=\rho$; when $\rho_c<\infty$, for all $\infty>\rho\ge \rho_c$, let $\varphi(\rho)=\phi_c$.

For any $j\in\N^*$, let
\[
F_j(\phi):=\sum_{r=j}^\infty rQ_r\phi^r,
\]
and
\[
\rho_{c,j}:=\sup_{0\le \phi<\phi_c}F_j(\phi)\in\overline{\R}_+.
\]
Note that the radius of convergence of series $F_j$ is $\phi_c$ and $F_j$ is smooth and strictly increasing on $(0,\phi_c)$ as well.
Similarly, for $x<\infty$ and $0\le x\le \rho_{c,j}$, let $\varphi_j(x)$ be the unique root of $F_j(\phi)=x$; for all $\infty>x\ge \rho_{c,j}$, let $\varphi_j(x)=\phi_{c}$.

\end{defi}

In this paper, we consider the case $\phi_c^{-1}:=\lim_{r\to\infty}Q_{r}^{-1/r}\in (0,+\infty)$. By the definition of $(Q_r)$, without loss any generality, we can always take $\phi_c>1$.

\begin{theorem}[Large deviation principle]\label{main:ldp}
Suppose $\frac{M}{V}\to \rho\in(0,\infty)$, then the scaled random variables
$$\bar{\eta}:=\frac{1}{V}\left(\eta_1,\dots,\eta_M\right)$$ under $\pi^{V,M}$ satisfies a large deviation principle on
$\R_+^\infty$
in the pointwise topology %weak~*~topology 
with speed $V$ and the good rate function
\[ J_\rho(\ell)= \left\{
\begin{array}{ll}
H(\ell|c^\rho)-\left(\rho-\sum_{r=1}^\infty r\ell_r\right)\log \frac{\varphi(\rho)}{\phi_c} & \text{if } \sum_{r=1}^\infty r \ell_r\le \rho,\\
\infty & \text{else},
\end{array} \right. \]
where   $c^\rho\in\R_+^\infty$ with $c_r^\rho=Q_r\varphi(\rho)^r$ for all $r\ge 1$ and  $H$ is the relative entropy function
  \[
H(\ell|c^\rho)=\sum_{r=1}^\infty\left(\ell_r\log\left(\frac{\ell_r}{Q_r\varphi(\rho)^r}\right)-\ell_r+Q_r\varphi(\rho)^r\right).
  \]

\end{theorem}

We now introduce a random variable $X$ whose density is totally determined by $(Q_r)$.
It will be shown later that the appearance of  one giant particle and the shape of it fluctuations depend on the tail behaviors of $X$. In other words, if the tail of the sequence $(Q_r)$ satisfies the following condition, when $M$ is large, there would be only one giant particle under the invariant measure $\pi^{V,M}$.

\begin{assumption}\label{basicapt}
  Suppose $q:=\sum_r Q_r\phi_c^r<\infty$ and $\rho_c=\sum_r rQ_r\phi_c^r<\infty$. Let $X$ be a random variable taking values in $\N^*:=\{1,2,\dots\}$, where
\[
\P(X=k)=\frac{Q_k\phi_c^k}{\sum_{r=1}^\infty Q_r\phi_c^r}.
\]

\begin{enumerate}
\item[(I)] Suppose the density of $X$ is subexponential
 \begin{equation}\label{sube}
\lim_{m\to\infty}\frac{\P(X_1+X_2=m)}{2\P(X=m)}=1.
\end{equation}
\item[(II)] Suppose there exist a sequence  $x(n)$, non-decreasing at infinity  such that
  \begin{equation}\label{asysub}
\lim_{n\to \infty}\sup_{m\ge  n\E X+x(n)}\left|\frac{\P(\sum_{i=1}^nX_i=m)}{n\P(X= \lceil m- n\E X\rceil)}-1\right|=0,
  \end{equation}
  where $X_i$ is a sequence of \emph{i.i.d} copies of $X$.
\item[(III)] Suppose the tail of $X$ is intermediate regularly varying,
\begin{equation}\label{inter}  \left\{
\begin{array}{l}
 \lim_{\eps\to 0^+}\liminf_{m\to\infty}\inf_{|k/m-1|<\eps}
  \frac{\P(X=k)}{\P(X=m)}=1,\\
 \lim_{\eps\to 0^+}\limsup_{m\to\infty}\sup_{|k/m-1|<\eps}
   \frac{\P(X=k)}{\P(X=m)}=1.
\end{array} \right.
\end{equation}
\end{enumerate}
\end{assumption}
\begin{remark}
For more about subexponential random variables, we refer to the book Foss et al.~\cite{foss}. The condition~\eqref{asysub} is related to the appearance of a giant jump in a random walk with heavy tail jumps. It has been shown by  Denisov et al.~\cite{denisov} that a subclass of subexponential random walk satisfies this condition. See also the rich references cited by Denisov et al.~\cite{denisov}.  For more about intermediate regularly varying random variables, we refer to the paper Cline~\cite{cline}.
\end{remark}

\begin{theorem}[Condensation]\label{main:cond}

  Let $K_V$ be the size of the largest particle,
  \[
K_V:=\max\left\{i\bigg|1\le i\le M,\eta_i\neq 0\right\}.
\]
Suppose the parameters $(Q_r)$ satisfy the Assumption~\ref{basicapt} and $M\ge V \rho_c +x(qV)$ for $V$ sufficient large,
then under $\pi^{V,M}$, as $V\to\infty$, we have, in probability,
\begin{equation}\label{lln}
  \lim_{V\to\infty}\frac{K_V-(M-\rho_cV)}{V}=0.
\end{equation}
Moreover, then for all fixed $k$, the vector
$$\left(\eta_1,\dots,\eta_k\right)$$ under $\pi^{V,M}$ satisfies a CLT limit,
   \begin{equation}\label{mgclt}
\frac{(\eta_i,1\le i\le k)-V(Q_i\phi_c^i,1\le i\le k)}{\sqrt{V}}\Rightarrow \cal{N}(0,\Sigma),
\end{equation}
for the convergence in distribution as $V\to\infty$, where $\cal{N}(0,\Sigma)$ is the $k$ dimensional Gaussian random variables with mean $0$ and variance $\Sigma={\rm diag}(Q_1\phi_c,Q_2\phi_c^2,\dots,Q_k\phi_c^k)$.
\end{theorem}
\begin{remark}
  If  $M/V\to \rho\in(\rho_c,\infty)$ and the parameters $(Q_r)$ satisfy the Assumption~\ref{basicapt}  with $\limsup_{n\to\infty}x(n)/n<\rho-\rho_c$, then we have, in probability,
  \[
\lim_{V\to\infty}\frac{K_V}{V}=\rho-\rho_c.  \]
\end{remark}

Under additional assumptions, we can obtain the central limit theorem for this condensed mass.
\begin{theorem}\label{main:clt}
  Suppose all the conditions in Theorem~\ref{main:cond} are satisfied.
\begin{enumerate}
\item[(a)] If in addition,
  \[
  \sum_{r=1}^\infty r^2Q_r\phi_c^r<\infty, \]
  then under $\pi^{V,M}$, as $V\to\infty$, we have, for the convergence in distribution, 
  \begin{equation}\label{gauss}
  \frac{K_V-(M-\rho_cV)}{\sqrt{V}}\Rightarrow \cal{N}\left(0,\sum_{r=1}^\infty r^2Q_r\phi_c^r\right).
  \end{equation}
\item[(b)]  Suppose there exists a slowly varying function $L:\N^*\mapsto \R_+$ and a constant $1<\alpha<2$, such that
\[
\sum_{j\ge r}Q_j\phi_c^j=\frac{L(r)}{r^{\alpha}}
\]
then under  $\pi^{V,M}$, as $V\to\infty$, we have, for the convergence in distribution, 
\begin{equation}\label{stable}
\frac{K_V-(M-\rho_cV)}{a_V}\Rightarrow \cal{L}_\alpha,
\end{equation}
where
\[
a_V:=\inf\left\{x:\P\left(\sum_{i=1}^{N(q)}X_i>x\right)\le \frac{1}{V}\right\}.
\]
and $\cal{L}_\alpha$ is a $\alpha-$stable random variable with characteristic function
\[
\E e^{it\cal{L}_\alpha}=\exp\left(itc-\int_0^\infty\left(e^{itx}-1-\frac{itx}{1+x^2}\right)\alpha x^{-(\alpha+1)}\diff x\right),
\]
where $c$ is a constant.
\item[(c)] Suppose $\sum_{r=1}^\infty r^2Q_r\phi_c^r=\infty$ and the law of $X$ belongs to the domain of attraction of a normal law,
  \[
\lim_{r\to\infty}\frac{r^2\sum_{j>r}Q_j\phi_c^j}{\sum_{i=1}^rj^2Q_j\phi_c^j}=0,
\]
then there exists a sequence $C_V\to\infty$ such that, under $\pi^{V,M}$, for the convergence in distribution,  we have
\begin{equation}\label{ga}
  \frac{K_V-(M-V\sum_{j=1}^{C_V} jQ_j\phi_c^j)}{\sqrt{V\sum_{j=1}^{C_V}j^2Q_j\phi_c^j}}\Rightarrow \cal{N}\left(0,1\right).
\end{equation}
\end{enumerate}
  \end{theorem}

\begin{example}[Power law case]
  From Nagaev~\cite{MR0282396}, Tka\v{c}uk~\cite{MR0368115} and Baltrunas~\cite{MR1429811}, we know that when the law of $X$ follows a Power law decay, that is, for some $b>2$,
  \[
\P(X=r)=\frac{ Q_r\phi_c^r}{\sum_jQ_j\phi_c^j}=\frac{1}{\zeta(b)}\frac{1}{r^b},
\]
where $\zeta(\cdot)$ is the Riemann Zeta function,
then the Assumption~\ref{basicapt} holds with $x(n)/n\to \eps$ for some $\eps>0$. See also  Denisov et al.~\cite{denisov} and its references. The power law condition has been considered in the Zero-Range model~\cite{armen, armendariz} as well as in the ideal Bose gas model~\cite{adams2007large,adams2008large}. An application of Theorem~\ref{main:cond} tells us that when $M/V\to \rho$ and $\rho>\zeta(b-1)+\eps$,  the size of the largest particle is  approximately $V(\rho-\zeta(b-1))$. Moreover, by using Theorem~\ref{main:cond}, we obtain the CLT for this largest particle. When $b>3$, the fluctuation is of order $\sqrt{V}$ and is approximately Gaussian with variance $\zeta(b-2)$. When $2<b<3$, the fluctuation is approximately a $(b-1)-$stable random variable of order $a_V=O(V^{1/(b-1)})$ by using Theorem~3.37 in Foss et al.~\cite{foss}.  When $b=3$,  the fluctuation is of order $\sqrt{V\log V}$ and is approximately Gaussian with variance $1/2$ by taking $C_V=V^{1/2}\log V$. See the proof of Theorem~\ref{main:clt} later for more details.

\end{example}

\section{Proof of the large deviation principle}\label{sec:ldp}
In order to prove the large deviation principle in Theorem~\ref{main:ldp}, we introduce a sequence of compound Poisson process that provides an alternative representation of the measure $\pi^{V,M}$. This approach can also be used to explain the appearance of the giant cluster in the super-critical sparse Erd\"os-R\'enyi graph. We refer to the paper by the author~\cite{sunmdp}.

\subsection{The invariant measure in terms of a conditional law}

Let $Y_1,Y_2,\dots$ be a sequence of \emph{i.i.d.} random variables on $N^*$ with common law
\[
\P(Y=r)=\frac{Q_r}{w},\qquad \forall r\in\N^*,
\]
where $w:=\sum Q_r<\infty$ by the assumption $\phi_c>1$.
Let $N(wV)$ be a Poisson process with intensity $wV$ that is independent of the sequence  $ (Y_i)$. For all $r\in\N^*$, let
\[
\cal{L}_r^V=\sum_{i=1}^{N(wV)}\ind{Y_i=r},
\]
be the $r-$th marginal empirical measure of our compound Poisson process $$\{Y_1,Y_2,\dots,Y_{N(wV)}\}.$$ We now show that the invariant measure $\pi^{V,M}$ can be expressed by a conditional law of this process. 

\begin{lemma}\label{lmcpp}
For all $\eta\in\N^M$ such that $\sum_{r=1}^Mr\eta_r=M$, we have
\[
\pi^{V,M}(\eta)=\P\left(\cal{L}_r^V=\eta_r,1\le \forall r\le M \bigg|\sum_{i=r}^{N(wV)}Y_r=M\right).
\]
\end{lemma}
\begin{proof}
  By the properties of Poisson point process, the random variables
 $(\cal{L}_r^V)_{1\le r\le M}$
are independent Poisson random variables with parameters $(VQ_r)_{1\le j\le M}$. Therefore,  for all  $\eta\in\N^M$ such that $\sum_{r=1}^Mr\eta_r=M$,
  \begin{equation*}
RHS=\frac{1}{\P\left(\sum_{r=1}^{N(wV)}Y_r=M\right)}\prod_{r=1}^M\frac{\left(VQ_r\right)^{\eta_r}}{\eta_r!}=\pi^{V,M}(\eta).
  \end{equation*}
\end{proof}

\subsection{Proof of Theorem~\ref{main:ldp}}

\begin{lemma}\label{lol}
  For any $j\in\N^*$, the sequence of average $\{\sum_{i=1}^{N(wV)}Y_i\ind{Y_i\ge j}/V\}$ satisfies a LDP in $\R_+$ with speed $V$ and the good rate function
  \[
I_{j}(x):
=x \log \left({\varphi_{j}(x)}\right)-\sum_{r=j}^\infty Q_r(\varphi_{j}(x))^r+\sum_{r=j}^\infty Q_r.
  \]
  Moreover, for any $\alpha>0$ and any sequence of integers $D\to\infty$ with $D/V\to\infty$, we have
  \begin{equation}\label{lmt}
\lim_{V\to\infty}\frac{1}{V}\log\P\left(\sum_{i=1}^{N(wV)}Y_i\ind{Y_i\ge j}=D\right)=-I_j(\alpha).
  \end{equation}

\end{lemma}
\begin{proof}

For all $k\ge 1$, the logarithmic moment generating function of the random variable $\sum_{i=1}^{N(wV)}Y_i\ind{Y_i\ge k}$ is given by
\[
\frac{1}{V}\log  \E\left(\exp\left(\theta \sum_{i=1}^{N(wV)}Y_i\ind{Y_i\ge k}\right)\right)
=\sum_{r=k}^\infty Q_r\left(e^{\theta r}-1\right),
\]
which is finite when $\theta<\log \phi_c$.
Since $\phi_c>1$, the large deviation principle follows from the classical Cram\'er's theorem with a good rate function,
\[
\sup_{\theta\in\R}\left\{x \theta -\sum_{r=k}^\infty Q_r e^{\theta r}\right\}+\sum_{r=k}^\infty Q_r
=x \log \left({\varphi_{k}(x)}\right)-\sum_{r=k}^\infty Q_r(\varphi_{k}(x))^r+\sum_{r=k}^\infty Q_r.
\]
To prove~\eqref{lmt}, firstly, by the Chebycheff inequality, we have the upper bound
\[
  \limsup_{V\to\infty}\frac{1}{V}\log\P\left(\sum_{i=1}^{N(wV)}Y_i\ind{Y_i\ge k}=D\right)\le\inf_{\theta\in\R}\left\{\sum_{r=k}^\infty Q_r(e^{\theta r}-1)-\theta\alpha\right\}=-I_k(\alpha).
  \]
  To prove the lower bound, for each $V$, we  define a $\zeta\in\N^{\N^*}$ such that $\sum_{r\ge j} r\zeta_r=D$.  For some $h\in(0,1/2)$ fixed, if
$k:=D-\sum_{r=1}^{\lfloor V^h\rfloor}r\lfloor V Q_r\varphi_j(D/V)^r\rfloor\notin [j, V^h]$, then let
  \[
\zeta_r:=\lfloor V Q_r\varphi_j(D/V)^r\rfloor,\qquad \forall j\le r\le V^{h},
\]
and $\zeta_k=1$. For all other index $r$,  let $\zeta_r=0$.
If $k\in [j,V^h]$, then let
\[
\zeta_r:=\lfloor V Q_r\varphi_j(D/V)^r\rfloor,\qquad \forall j\le r\neq k\le V^{h},
\]
\[
\zeta_k:=\lfloor V Q_r\varphi_j(D/V)^k\rfloor+1,
\]
and $\zeta_r=0$ for all $r\notin [j,V^h]$.

Let $\theta:=\sum_{r=j}^\infty \zeta_r$, then we have
\begin{multline*}
\P\left(\sum_{i=1}^{N(wV)}Y_i\ind{Y_i\ge j}=D\right)\\
\ge \sum_{n\ge \theta}\P\left(N(wV)=n\right)\frac{n!}{(n-\theta)!}\left(\P(X<j)^{n-\theta}\right)\prod_{r=j}^\infty\frac{1}{\zeta_r!}\P(X=r)^{\zeta_r}\\
=e^{-wV\P(X\ge j)}\prod_{r=j}^{\infty}\frac{(VQ_r)^{\zeta_r}}{\zeta_r!}.
\end{multline*}
By using the facts $k/V\to (\alpha-\rho_{c,j})^+$ and $\lim_{k\to\infty}Q_k^{-1/k}=\phi_c$, we have
\[
\lim_{V\to \infty}\frac{1}{V}\log Q_k=-(\alpha-\rho_{c,j})^+\log \phi_c.
\]
Combining with the Stirling's formula, we obtain
\begin{multline*}
\liminf_{V\to\infty}\frac{1}{V}\log\P\left(\sum_{i=1}^{N(wV)}Y_i\ind{Y_i\ge j}=D\right)\\
\ge-\sum_{r=j}^\infty Q_r+\sum_{r=j}^{\infty} Q_r\varphi(\alpha)^r\left(1-r\log \varphi(\alpha)\right)-(\alpha-\rho_{c,j})^+\log \phi_c\\
=-\sum_{r=j}^\infty Q_r+\sum_{r=j}^{\infty} Q_r\varphi(\alpha)^r- \alpha\log( \varphi(\alpha))=-I_j(\alpha).
\end{multline*}
In conclusion, we have
\[
\lim_{V\to\infty}\frac{1}{V}\log\P\left(\sum_{i=1}^{N(wV)}Y_i\ind{Y_i\ge j}=D\right)=-I_j(\alpha).
\]

\end{proof}
\begin{remark}
When $j=1$, $I_1(\rho)$ is related to the free energy of the measure $\pi^{V,M}$.
\end{remark}

\begin{prop}
  [Finite marginal LDP]\label{finiteldp}
Fix $j\in \N^*$, when $\frac{M}{V}\to \rho\in(0,\infty)$, the scaled random variables
$$\frac{1}{V}\left(\cal{L}^V_1,\cal{L}^V_2,\dots,\cal{L}^V_j\right)$$ under condition $\sum_{i=1}^{N(wV)}Y_i=M$  satisfies a large deviation principle on $\R_+^j$
with speed $t$ and the good rate function
\begin{multline*}
  J_{\rho,j}(\ell_1,\dots,\ell_j)\\
  = \sum_{r=1}^j\ell_r\left(\log \frac{\ell_r}{Q_r}-1\right)+(\rho-m_j) \log \left({\varphi_{j+1}(\rho-m_j)}\right)-\sum_{r=j+1}^\infty Q_r(\varphi_{j+1}(\rho-m_j))^r\\
  -\rho \log \left({\varphi(\rho)}\right)+\sum_{r=1}^\infty Q_r(\varphi(\rho))^r,
\end{multline*}
 if $m_j:=\sum_{r=1}^jr\ell_r\le \rho$ and $J_{\rho,j}(\ell_1,\dots,\ell_j)=\infty$ if else.

\end{prop}

\begin{proof}

  For any $\eps>0$, there exists $V$ sufficient large such that
  \[
\P\left(\frac{1}{V}(\cal{L}^V_r,1\le r\le j)\in \cal{B}_{\rho+\eps}^j\bigg|\sum_{i=1}^{N(wV)}Y_i=M\right)=1,
\]
where
\[
\cal{B}_{\rho+\eps}^j:=\left\{\ell\in\R_+^j\bigg|\sum_{r=1}^jr\ell_r\le \rho+\eps\right\},
\]
which is a compact set in $\R_+^j$. Hence, we only need to prove the weak large deviation principle.
  
For any $\ell=(\ell_r,1\le r\le j)\in \R_+^j$  and any $\delta>0$, let ${B_\delta(\ell)}$(\emph{resp.} $\overline{B_\delta(\ell)}$) be the open (\emph{resp.} closed) ball centered at $\ell$ with radius $\delta$ in $\R_+^j$. If $\sum_{r=1}^jr\ell_r>\rho$, then for any $\eps\in (0,\sum_{r=1}^jr\ell_r-\rho)$, there exists $\delta$ small enough, such that $\inf_{c\in\overline{ B_\delta(\ell)}}\sum_{r=1}^jrc_r\ge \rho+\eps$. Therefore, for $V$ sufficient large, we have
\[
 \P\left(\frac{1}{V}(\cal{L}^V_r,1\le r\le j)\in \overline{B_\delta(\ell)}\bigg|\sum_{i=1}^{N(wV)}Y_i=M\right)=0.
 \]
 If $\sum_{r=1}^jr\ell_r\le \rho$, the limit
 \[
 \lim_{\delta\to 0}\limsup_{V\to\infty} \frac{1}{V}\log \P\left(\frac{1}{V}(\cal{L}^V_r,1\le r\le j)\in \overline{B_\delta(\ell)}\bigg|\sum_{i=1}^{N(wV)}Y_i=M\right)
 \]
 can be upper bounded by
 \begin{multline*}
   \lim_{\delta\to 0}\limsup_{V\to\infty}  \frac{1}{V}\log \P\left(\frac{1}{V}(\cal{L}^V_r,1\le r\le j)\in \overline{B_\delta(\ell)}\right)\hfill\\
   +\lim_{\delta\to 0}\limsup_{V\to\infty} \frac{1}{V}\log \P\left(\frac{1}{V}\sum_{i=1}^{N(wV)}Y_i\ind{Y_i>j}- \frac{M}{V}+\sum_{r=1}^j r\ell_r\in \overline{B}_{j^{3/2}\delta}(0)\right)
   \\
   -\lim_{V\to\infty}\frac{1}{V}\log \P\left(\sum_{i=1}^{N(wV)}Y_i=M\right).\\
\end{multline*}
Recall that $\cal{L}^V_r$, $1\le r\le j$, are independent Poisson random variables with parameter $VQ_r$. Therefore,  by using the large deviation principles of the scaled Poisson random vector$(\cal{L}^V_r/V, 1\le r\le j)$  and Lemma~\ref{lol},  this upper bound equals to
\[
   -\sum_{r=1}^jH(\ell_r|Q_r)-I_{j+1}(\rho-m_j)+I_1(\rho)
   =-J_{\rho,j}(\ell_1,\dots,\ell_j),
   \]
   where
   \[
H(\ell_r|Q_r)=\ell_r\log\frac{\ell_r}{Q_r}-\ell_r+Q_r.
   \]
On the other hand, for the similar reason, the limit
\[
 \lim_{\delta\to 0}\liminf_{V\to\infty} \frac{1}{V}\log \P\left(\frac{1}{V}(\cal{L}^V_r,1\le r\le j)\in {B_\delta(\ell)}\bigg|\sum_{i=1}^{N(wV)}Y_i=M\right)
\]
is lower bounded by
\begin{multline*}
 \liminf_{V\to\infty} \frac{1}{V}\log \P\left(\cal{L}^V_r=\lfloor V\ell_r\rfloor,1\le r\le j\right)\hfill\\+\lim_{V\to\infty}\frac{1}{V}\log \P\left(\sum_{i=1}^{N(wV)}Y_i\ind{Y_i>j}=M-\sum_{r=1}^j r\lfloor V\ell_r\rfloor\right)
   \\ -\lim_{V\to\infty}\frac{1}{V}\log \P\left(\sum_{i=1}^{N(wV)}Y_i=M\right)=-J_{\rho,j}(\ell_1,\dots,\ell_j).
  \end{multline*}
Finally, the LDP holds by applying   the G\"artner-Ellis theorem.
To show that the rate function is good, we firstly notice that the relative entropy function $\sum_{r=1}^jH(\ell_r|Q_r)$ is good on $\R_+^j$. From Lemma~\ref{lol}, we find $I_{j+1}(\cdot)$ is a good rate function on $\R_+$. Since the mapping $(\ell_r,1\le r\le j)\mapsto \rho-\sum_{r=1}^jr\ell_r$ is continuous on $\R^j$, the rate function $I_{j+1}(\rho-\sum_{r=1}^j r\ell_r)$ is also good on $\R^j_+$. Therefore, for any $\alpha>0$, the level set $\{\ell\in\R_+^j|J_{\rho,j}(\ell_1,\dots,\ell_j)\le \alpha\}$ is a closed set in $\R_+^j$. Moreover, by the relation
\[
\left\{\ell \in\R_+^j\bigg|J_{\rho,j}(\ell_1,\dots,\ell_j)\le \alpha\right\}\subset \left\{\ell \in\R_+^j\bigg|\sum_{r=1}^jH(\ell_r|Q_r)\le \alpha+I_1(\rho)\right\},
\]
this level set is bounded and hence a compact set in $\R^j$. In conclusion, $J_{\rho,j}(\ell_1,\dots,\ell_j)$ is a good rate function. 

\end{proof}

\begin{prop}\label{ratefun}
  For all $\ell\in\R_+^\infty$, we have

\[ \sup_{j\ge 1} J_{\rho,j}(\ell_1,\dots,\ell_j)= \left\{
\begin{array}{ll}
H(\ell|c^\rho)-\left(\rho-\sum_{r=1}^\infty r\ell_r\right)\log \frac{\varphi(\rho)}{\phi_c} & \text{if } \sum_{r=1}^\infty r \ell_r\le \rho,\\
\infty & \text{else}.
\end{array} \right. \]
\end{prop}

\begin{proof}
For all $\ell\in\R_+^\infty$ with $m:=\sum_{r=1}^{\infty}r\ell_r\le \rho$, $m_j=\sum_{r=1}^jr\ell_r$, the difference,
  \[
    \Delta_j:=H(\ell|c^\rho)+\left(\sum_{r=1}^\infty r \ell_r-\rho\right)\log\frac{\varphi(\rho)}{\phi_c}-J_{\rho,j}(\ell_1,\dots,\ell_j),\]
equals to 
\begin{multline*}
    \sum_{r=j+1}^\infty\left(\ell_r\log\left(\frac{\ell_r}{Q_r\varphi_{j+1}(\rho-m_j)^r} \right)  -\ell_r+Q_r\left(\varphi_{j+1}(\rho-m_j)\right)^r\right) \\
      -(\rho-m)\log\left(\frac{\varphi_{j+1}(\rho-m_j)}{\phi_c}\right),\\
  \end{multline*}
 which  is non-negative. Therefore, we have
  \[
 H(\ell|c^\rho)+\left(\sum_{r=1}^\infty r \ell_r-\rho\right)\log\frac{\varphi(\rho)}{\phi_c}\ge \sup_{j\ge 1}J_{\rho,j}(\ell_1,\dots,\ell_j).
  \]
  On the other hand,
  by simple calculation, the difference $\Delta_j$ also equals to
 \begin{multline}\label{delj}
    \sum_{r=j+1}^\infty\left(\ell_r\log\left(\frac{\ell_r}{Q_r}\right)    -\ell_r\right)
      -(\rho-m_j)\log\left(\frac{\varphi_{j+1}(\rho-m_j)}{\phi_c}\right)\\+\sum_{r=j+1}^\infty Q_r\left(\varphi_{j+1}(\rho-m_j)\right)^r-(m-m_j)\log \phi_c.
 \end{multline}
 We now show that when $H(\ell|c^\rho)<\infty$, for any $\eps>0$ and for $j$ sufficient large, the difference $\Delta_j<\eps$. Firstly, it is easy to see that the facts $H(\ell|c^\rho)<\infty$, $m=\sum_r r\ell_r\le \rho<\infty$, $\varphi(\rho)\le \phi_c<\infty$ and $\sum_rQ_r<\infty$ imply
 \[
\sum_{r=1}^\infty\left(\ell_r\log\left(\frac{\ell_r}{Q_r}\right)    -\ell_r+Q_r\right)<\infty.
\]
Hence the first term in \eqref{delj} is vanishing when $j\to\infty$. If $\rho-m>0$, then the sequence $\varphi_{j+1}(\rho-m)$ is increasing. Moreover, 
if $\lim_{j\to\infty}\varphi_{j+1}(\rho-m)<\phi_c$, then
\[
\rho-m=F_{j+1}(\varphi_{j+1}(\rho-m))\le F_{j+1}(\lim_{j\to\infty}\varphi_{j+1}(\rho-m))\to 0,
\]
that leads to a contradiction. Thus, by the relation
\[
\phi_c\ge \varphi_{j+1}(\rho-m_j)\ge \varphi_{j+1}(\rho-m),
\]
we must have
\[
\lim_{j\to\infty}\varphi_{j+1}(\rho-m_j)=\phi_c.
\]
Hence the second term in \eqref{delj} is vanishing.
If $\rho-m=0$, by using
  \[
  \phi_c\ge \varphi_{j+1}(\rho-m_j)\ge \varphi_{1}(\rho-m_j)
  \]
  and $\varphi_{1}(\rho-m_j)\to 0$, we find that, for $j$ sufficient large, there exists a $\phi'<\phi_c$ such that  $\varphi_{1}(\rho-m_j)<\phi'$ and
  \[
F_1(\phi')\frac{\varphi_{1}(\rho-m_j)}{\phi'}\ge F_1(\varphi_{1}(\rho-m_j))=(\rho-m_j).
\]
Thus,
\[
(\rho-m_j)\log\left(\varphi_{1}(\rho-m_j)\right)\ge (\rho-m_j)\log\left(\frac{\phi'}{F_1(\phi')}(\rho-m_j)\right).
\]
Therefore,
\[
\liminf_{j\to\infty}(\rho-m_j)\log   \varphi_{j+1}(\rho-m_j)\ge
\lim_{j\to \infty}(\rho-m_j)\log\left(\frac{\phi'}{F_1(\phi')}(\rho-m_j)\right)=0.
\]
On the other hand, it is smaller than
\[
\limsup_{j\to\infty}(\rho-m_j)\log   \varphi_{j+1}(\rho-m_j)\le  \lim_{j\to\infty}(\rho-m_j)\log   \phi_c=0.
\]
Hence, the second term in~\eqref{delj} is vanishing.  The third term is vanishing since
\begin{multline*}
\sum_{r=j+1}^\infty Q_r\left(\varphi_{j+1}(\rho-m_j)\right)^r\le\frac{1}{j+1} \sum_{r=j+1}^\infty rQ_r\left(\varphi_{j+1}(\rho-m_j)\right)^r\\
\le \frac{1}{j+1} \left((\rho-m_j)\land \rho_{c,j}\right)\to 0 \qquad as~j\to\infty.
\end{multline*}
Combining with the fact $m_j\to m$, we obtain
\[
H(\ell|c^\rho)+\left(\sum_{r=1}^\infty r \ell_r-\rho\right)\log\frac{\varphi(\rho)}{\phi_c}\le \sup_{j\ge 1}J_{\rho,j}(\ell_1,\dots,\ell_j)+\eps,
\]
for any $\eps>0$. We conclude the proof by letting $\eps\to 0$.

If $H(\ell|c^\rho)=\infty$, then
\[
\sup_{j\ge 1} J_{\rho,j}(\ell_1,\dots,\ell_j)\ge
\sum_{r=1}^j\ell_r\left(\log \frac{\ell_r}{Q_r}-1\right)+(\rho-m_j) \log \left({\varphi_{j+1}(\rho-m_j)}\right)
  +\sum_{r=1}^\infty Q_r(\varphi(\rho))^r,
  \]
  for all $j$. From previous proofs, we see that the last two terms are finite. Therefore, must we have
  \[\sup_{j\ge 1} J_{\rho,j}(\ell_1,\dots,\ell_j)=\infty.\]
  When $\sum_{r=1}^\infty r\ell_r>\rho$, then there exists a finite $k$ such that $\sum_{r=1}^k r\ell_r>\rho$. In this case,
  \[
  \sup_{j\ge 1} J_{\rho,j}(\ell_1,\dots,\ell_j)\ge
 J_{\rho,k}(\ell_1,\dots,\ell_k)=\infty.
  \]

\end{proof}

\begin{proof}[Proof of Theorem~\ref{main:ldp}]
 By using Proposition~\ref{finiteldp}, Proposition~\ref{ratefun} and applying Dawson-G\"artner Theorem (Theorem~4.6.1 in Dembo and Zeitouni~\cite{dembo}), we obtain the result.
\end{proof}

\section{On the giant jump and the equivalence of ensembles in a compound Poisson process}\label{sec:cond}
In this section, we introduce a new sequence of compound Poisson process in order to study the condensation phenomenon later.
Let $\{X_i\}_{i\in\N}$ be a sequence of \emph{i.i.d.} random variables on $\N^*$ with law
\[
\P(X=r)=\frac{Q_r\phi_c^r}{q},
\]
where $q=\sum rQ_r\phi_c^r<\infty$. Recall that under the Assumption~\ref{basicapt}, $\rho_c=q\E X<\infty$. We remark that in the previous section, $\rho_c$ can be infinite. Let $N(qV)$ be a Poisson process with intensity $qV$ and be independent of the sequence $(X_i)$. By a similar proof as  Lemma~\ref{lmcpp}, we can easily get another representation of our invariant measure as follows.

\begin{lemma}\label{lmcpp2}
For all $\eta\in\N^M$ such that $\sum_{r=1}^Mr\eta_r=M$, we have
\[
\pi^{V,M}(\eta)=\P\left(\sum_{i=1}^{N(qV)}\ind{X_i=r}=\eta_r,1\le \forall r\le M \bigg|\sum_{i=1}^{N(qV)}X_r=M\right).
\]
\end{lemma}

We will show that, under the Assumption~\ref{basicapt} and the condition
\[
\sum_{i=1}^{N(qV)}X_r=M,
\]
with $M$ sufficient large, when $V\to\infty$, there is only one giant jump appears in the compound Poisson $\{X_1,X_2,\dots,X_{N(qV)}\}$. Moreover, the rest of the jumps are asymptotically \emph{i.i.d.} distributed. Similar phenomenon happens in the random walk under a subclass of subexponential condition. See the papers by Denisov et al.~\cite{denisov} and Armend\'ariz et al.~\cite{armen}.

To describe the jumps but the largest one, we need to define a mapping
\[
{\rm bulk}: \cup_{n\ge 1}\N^n\mapsto \cup_{n\ge 0}\N^n,
\]
with convention $N^0=\emptyset$, by
\[
{\rm bulk}(x_1,x_2,\dots,x_n)=(x_1,\dots,x_{m*-1},x_{m*+1},\dots,x_n),
\]
where
\[
m^*:=\inf_{i\ge 1}\{x_i=\max_{1\le j\le n}x_j\}.
\]

In the rest of the section, we will prove the following results.

\begin{prop}\label{poisum}
  Suppose the law of $X$ satisfies the Assumption~\ref{basicapt},
  then we have 
\begin{equation}\label{thm:ld}
\lim_{V\to\infty}\sup_{m\ge y(qV)} \left|\frac{\P\left(\sum_{i=1}^{N(qV)}X_i=m\right)}{ qV\P\left(X=\lceil m-\rho_c V\rceil\right)}-1\right|=0,
\end{equation}
where
\[ y(t):= \left\{
\begin{array}{ll}
\gamma t & \text{if } \xi<\infty,\\
(1+\omega)x((1+\omega)t) & \text{else,}
\end{array} \right. \]
with $\xi:=\limsup_{n\to\infty}x(n)/n\in[0,+\infty]$,  $\gamma>(\rho_c/q+\xi)$ and $\omega>0$.

Moreover, we have
  \begin{multline}\label{thm:cond}
    \lim_{V\to\infty}\sup_{m\ge y(qV)}\sup_{A\in\cal{B}(\cup_{n\ge 1}(\N^*)^n)}\bigg|\P\left({\rm{bulk}}(X_1,\dots,X_{N(qV)})\in A\bigg|\sum_{i=1}^{N(qV)}X_i=m\right)\\-\P\left((X_1,\dots,X_{N(qV)-1})\in A\right)\bigg|= 0.
\end{multline}

\end{prop}

\subsection{Proof of Proposition~\ref{poisum}}
The proof of relation~\eqref{thm:ld} follows a similar approach in \cite{kumi} by showing the next four lemmas.

\begin{lemma}\label{int}
  Suppose the tail density of $X$ is intermediate regularly varying, cf. condition~\eqref{inter}, then we have,
   \begin{equation}\label{int1}
\limsup_{V\to\infty}\sup_{m\ge y(qV)}\frac{\P(X=m)}{\P(X=\lceil m-\rho_c V \rceil)}<\infty.
  \end{equation}
Moreover, for any $\eps\in(0,1)$, there exists $\delta_0\in(0,\eps/3)$, such that for $V$ sufficient large, we have
  \begin{equation}\label{int2}
\limsup_{V\to\infty}\sup_{m\ge y(qV)}\sup_{ |k-  qV|\le \delta_0 qV}\left|\frac{k\P(X=\lceil m-k\rho_c/q \rceil)}{  qV \P(X=\lceil m-\rho_c V\rceil)}-1\right|<\eps.
\end{equation}

\end{lemma}
\begin{proof}

By the intermediate regularly varying property~\eqref{inter}, for any $\eps>0$, there exists $\eps_0$, such that,
  \begin{equation}\label{ep3}
  \limsup_{s\to\infty}\sup_{|r/s-1|\le\eps_0}\left|\frac{\P(X=r)}{\P(X=s)}-1\right|<\eps/3.
  \end{equation}
  When $y(t)\gg t$, the relation~\eqref{int1} holds since
  \[
\frac{m}{m-\rho_cV}\to 1,\qquad \forall m\ge y(qV).
  \]
  For the case $y(t)=\gamma t$,
  we have
  \[
1\le \inf_{m\ge y(qV)}\frac{m}{\lceil m-k\rho_c/q \rceil}\le \sup_{m\ge y(qV)}\frac{m}{\lceil m-k\rho_c/q \rceil}\le \frac{\gamma}{\gamma-\rho_c/q}.
  \]
  Therefore,
  \[
  \limsup_{V\to\infty}\sup_{m\ge y(qV)}\frac{\P(X=m)}{\P(X=\lceil m-\rho_c V \rceil)}
  \le \limsup_{s\to\infty}\sup_{1\le r/s\le  \frac{\gamma}{\gamma-\rho_c/q}}\frac{\P(X=r)}{\P(X=s)},
  \]
  which is finite by iterating~\eqref{ep3}.

To prove~\eqref{int2}, we first notice that there exists $\delta_0\in(0,\eps/3)$, such that for $V$ sufficient large
\[
\sup_{m\ge y(qV)}\sup_{|k-  qV|\le \delta_0 qV}\left|\frac{\lceil m-k\rho_c/q\rceil}{\lceil m-\rho_c V\rceil}-1\right|\le \frac{\delta_0 \rho_c V+1}{m-\rho_c V}\le \frac{\delta_0 \rho_c V+1}{y(qV)-\rho_c V}<\eps_0.
\]
Hence, we have that
\begin{multline*}
\limsup_{V\to\infty}  \sup_{m\ge y(qV)}\sup_{|k-  qV|<\delta_0 qV}\left|\frac{k\P(X= \lceil m- k\rho_c/q\rceil))}{qV\P(X=\lceil m-\rho_c V\rceil)}-1\right|\\
  \le (1+\delta_0)\limsup_{s\to\infty}\sup_{|r/s-1|<\eps_0}\left|\frac{\P(X=r)}{\P(X=s)}-1\right|+\delta_0<\eps.
\end{multline*}

\end{proof}

\begin{lemma}\label{ldpdelta}
  For any sufficient small $\delta_0>0$, we have
  \[
\lim_{V\to 0}\sup_{m\ge y(qV)} \sup_{ |k-  qV|\le \delta_0 qV}\left|
\frac{\P\left(\sum_{i=1}^{k}X_i=m\right)}{k\P(X=\lceil m-k\rho_c/q \rceil)}-1\right|=0.
\]
\end{lemma}
\begin{proof}
If $\xi=\limsup x(n)/n<\infty$, then for any $\delta_0<(\gamma/( \rho_c/q+\xi)-1 )$,  for all $V$ sufficient large and all $k$ such that $|k-qV|\le \delta_0 qV$, we have
\[
k\rho_c/q +x(k)\le \gamma t=y(qV).
\]
If $x(n)\gg n$, then for all $\delta_0<\omega$, $V$ sufficient large and $|k-qV|\le \delta_0 qV$, we have
\[ k\rho_c/q+x(k)\le (1+\omega)x((1+\omega)qV)=y(qV).\]
In both cases, we have
\begin{multline*}
  \sup_{m\ge y(qV)}\sup_{|k-qV|\le\delta_0 qV}\left|\frac{\P\left(\sum_{i=1}^{k}X_i=m\right)}{k\P(X= \lceil m- k\rho_c/q\rceil))}-1\right|\\
  \le \sup_{k\ge (1-\delta_0)qV}  \sup_{m\ge   k\rho_c/q+x(k)}\left|\frac{\P\left(\sum_{i=1}^{k}X_i=m\right)}{k\P(X= \lceil m- k\rho_c/q\rceil))}-1\right|\to 0,\qquad as~V\to \infty,
\end{multline*}
thanks to the assumption~\eqref{asysub}.

\end{proof}

\begin{lemma}\label{lgdelta}

  For any $\eps>0$, there exists a small enough $\delta_0>0$, such that
\[
  \limsup_{V\to\infty}\sup_{m\ge y(qV)}\left|  \sum_{ |k-qV|<\delta_0 qV}
  \left(\frac{\P\left(\sum_{i=1}^{k}X_i=m,N(qV)=k\right)}{  qV \P(X=\lceil m-\rho_c V\rceil)}- \P(N(qV)=k)\right)\right|<\eps.
  \]

\end{lemma}

\begin{proof}
  By using Lemma~\ref{int} and Lemma~\ref{ldpdelta}, we have
  \begin{multline*}
 \limsup_{V\to\infty}\sup_{m\ge y(qV)}\left|  \sum_{ |k- qV|<\delta_0qV}
 \left(\frac{\P\left(\sum_{i=1}^{k}X_i=m,N(qV)=k\right)}{ qV \P(X=\lceil m-\rho_c V\rceil)}- \P(N(qV)=k)\right)\right|\\
 \le \limsup_{V\to\infty}\sup_{m\ge y(qV)} \sup_{ |k- qV|<\delta_0qV}\left| 
 \left(\frac{\P\left(\sum_{i=1}^{k}X_i=m\right)}{k\P(X=\lceil m- k\rho_c/q\rceil)}\frac{k\P(X=\lceil m-  k\rho_c/q\rceil)}{ qV \P(X=\lceil m-\rho_c V\rceil)}- 1\right)\right|<\eps.
    \end{multline*}
\end{proof}

\begin{lemma}\label{de0}

For any  $\delta_0\in(0,1 )$,  we have
  \[
\lim_{V\to \infty}\sup_{m\ge y(qV)}\left|  \sum_{ |k- qV|\ge \delta_0qV}\frac{\P\left(\sum_{i=1}^{k}X_i=m,N(qV)=k\right)}{ qV \P(X=\lceil m-\rho_c V\rceil)}\right|=0.
  \]
\end{lemma}
\begin{proof}
When the density of $X$ is subexponential, the Kesten's bound holds, \emph{i.e.}, for any $\eps>0$, there exists $c(\eps)>0$, such that for all $m,k\in\N$,
\[
\P\left(\sum_{i=1}^k X_i=m\right)\le c(\eps)(1+\eps)^k\P(X=m).
\]
Hence, we have upper bound
\begin{multline*}
\left|  \sum_{ |k- qV|\ge \delta_0qV}\frac{\P\left(\sum_{i=1}^{k}X_i=m,N(qV)=k\right)}{ qV \P(X=\lceil m-\rho_c V\rceil)}\right|\\
\le c(\eps)\frac{\E\left((1+\eps)^{N(qV)}\ind{|N(qV)- qV|\ge \delta_0qV}\right)}{ qV}\frac{\P(X=m)}{\P(X=\lceil m-\rho_c V\rceil)}.
\end{multline*}

For the term concerning the Poisson process, it is easy to see that,
\begin{multline*}
  \E\left((1+\eps)^{N(qV)}\ind{|N(qV)- qV|\ge \delta_0qV}\right)\\=
  \E\left((1+\eps)^{N(qV)}\ind{N(qV)\ge  qV + \delta_0qV}\right)+  \E\left((1+\eps)^{N(qV)}\ind{N(qV)\le  qV - \delta_0qV}\right),
\end{multline*}
which is smaller than
\begin{multline*}
  e^{-s( qV+\delta_0 qV)}\E\exp\{(s+\ln(1+\eps) )N(qV)\}+   e^{r( qV-\delta_0 qV)}\E\exp\{(\ln(1+\eps)-r)N(qV)\}\\
  =\exp\{ qV((1+\eps)e^s-1)-s( qV+\delta_0 qV)\}+\exp\{r( qV-\delta_0 qV)+ qV((1+\eps)e^{-r}-1)\},
\end{multline*}
 for any $s,r>0$.
For any $\eps<\delta_0<1$, let $s=\ln(\frac{1+\delta_0}{1+\eps})$ and  $r=\ln(\frac{1+\eps}{1-\delta_0})$, then we have
\[
\E\left((1+\eps)^{N(qV)}\ind{|N(qV)- qV|\ge \delta_0qV}\right)
\le \exp\left\{f_+(\eps,\delta_0)qV\right\}+\exp\left\{f_-(\eps,\delta_0)
qV\right\},
\]
where
\[
f_{\pm}(\eps,x)=-(1\pm x)\ln\left(\frac{1\pm x}{ 1+\eps}\right)\pm x.
\]
Clearly, for $\eps=0$, we have $f_+(0,\delta_0)<0$ and $f_-(0,\delta_0)<0$. Thus, there exists $\eps$ sufficient small such that $f_+(\eps,\delta_0)<0$ and $f_-(\eps,\delta_0)<0$. For such $\eps>0$, we have
\[
\lim_{V\to\infty}\sup_{m\ge y(qV)}c(\eps)\frac{\E\left((1+\eps)^{N(qV)}\ind{|N(qV)- qV|\ge \delta_0t}\right)}{ qV}\frac{\P(X=m)}{\P(X=\lceil m-\rho_c V\rceil)}=0,
\]
thanks to the relation~\eqref{int1}. 
\end{proof}

We now give the proof of Proposition~\ref{poisum}.
\begin{proof}[Proof of the limit~\eqref{thm:ld}]
  Since
    \begin{multline*}
 \left|\frac{\P\left(\sum_{i=1}^{N(qV)}X_i=m\right)}{ qV\P\left(X=\lceil m-\rho_c V\rceil\right)}-1\right|\\
    \le \left|  \sum_{ |k-qV|<\delta_0 qV}
    \left(\frac{\P\left(\sum_{i=1}^{k}X_i=m,N(qV)=k\right)}{  qV \P(X=\lceil m-\rho_c V\rceil)}- \P(N(qV)=k)\right)\right|\\
 +  \left|  \sum_{ |k- qV|\ge \delta_0qV}\frac{\P\left(\sum_{i=1}^{k}X_i=m,N(qV)=k\right)}{ qV \P(X=\lceil m-\rho_c V\rceil)}\right|
    +\P\left(|N(qV)-qV|\ge \delta_0 qV\right),
  \end{multline*}
then by using Lemma~\ref{lgdelta} and Lemma~\ref{de0}, for any $\eps>0$,  there exists $\delta_0\in(0,1)$, such that
\[
    \limsup_{V\to\infty}\sup_{m\ge y(qV)} \left|\frac{\P\left(\sum_{i=1}^{N(qV)}X_i=m\right)}{ qV\P\left(X=\lceil m-\rho_c V\rceil\right)}-1\right|
    <\eps.
 \]
  The limit~\eqref{thm:ld} holds by taking $\eps\to 0$.
\end{proof}

%%%%%%%%%%%%%%%%%%%%%%%%%%%%%%%%%
%%%%%%%%%%%%%%%%%%%%%%%%%%%%%%%%%
To show the limit ~\ref{thm:cond}, we will need the equivalence of ensemble result for the random walk $(X_1,\dots,X_n)$ from the paper  Armend\'ariz et al.~\cite{armen}.
\begin{lemma}\label{cutcut}
Suppose the law of $X$ satisfies the conditions~\eqref{sube} and~\eqref{asysub},
then we have
 \begin{multline*}
   \lim_{n\to\infty}\sup_{m\ge  n\E X  + x(n)}\sup_{A\in\cal{B}((\N^*)^{n-1}}\bigg|\P\left({\rm{bulk}}(X_1,\dots,X_{n})\in A\bigg|\sum_{i=1}^{n}X_i=m\right)\\
   -\P\left((X_1,\dots,X_{n-1})\in A\right)\bigg|= 0.
 \end{multline*}
\end{lemma}

\begin{proof}[Proof of the limit~\eqref{thm:cond}]
  By using Lemma~\ref{int}, Lemma~\ref{ldpdelta} and the limit~\eqref{thm:ld}, we have that, for any $\eps>0$, there exists $\delta_0\in(0,\eps/3)$, such that
\[
\limsup_{V\to\infty}\sup_{m\ge y(qV)}\sup_{|k-qV|\le\delta_0qV}\left|\frac{\P\left(\sum_{i=1}^{k}X_i=m\right)}{\P\left(\sum_{i=1}^{N(qV)}X_i=m\right)}-1\right|<2\eps.
\]
By the definition of $y(qV)$, we can choose $\delta_0$ small enough and $V$ large enough such that $y(qV)> k\rho_c /q+x(k)$ for all $|k-qV|<\delta_0qV$. Therefore, for $V$ large enough,
  \begin{multline*}
 \sup_{m\ge y(qV)} \sup_{A\in\cal{B}(\cup_{n\ge 1}(\N^*)^n)} \bigg|\P\left({\rm{bulk}}(X_1,\dots,X_{N(qV)})\in A,\left|N(qV)-qV\right|<\delta_0qV\bigg|\sum_{i=1}^{N(qV)}X_i=m\right)\\-\P\left((X_1,\dots,X_{N(qV)-1})\in A,\left|N(qV)-qV\right|< \delta_0qV\right)\bigg|\\
    \le   \sup_{m\ge y(qV)}  \sum_{\left|k-qV\right|< \delta_0qV}\sup_{A\in\cal{B}((\N^*)^{k-1})}\\\bigg|\frac{\P\left({\rm{bulk}}(X_1,\dots,X_{k})\in A,\sum_{i=1}^{k}X_i=m\right)}{\P\left(\sum_{i=1}^{N(qV)}X_i=m\right)}
    -\P\left((X_1,\dots,X_{k-1})\in A\right)\bigg|\P(N(qV)=k)\\
\le \sup_{k\ge (1-\delta_0)qV}\sup_{m\ge k\rho_c /q+x(k)}\sup_{A\in\cal{B}((\N^*)^{k-1})}\\\bigg|\frac{\P\left({\rm{bulk}}(X_1,\dots,X_{k})\in A,\sum_{i=1}^{k}X_i=m\right)}{\P\left(\sum_{i=1}^{k}X_i=m\right)}
-\P\left((X_1,\dots,X_{k-1})\in A\right)\bigg|(1+2\eps)+2\eps,
  \end{multline*}
which is upper bounded by $2\eps$ for $V\to \infty$ by using Lemma~\ref{cutcut}

On the other hand, thanks to Lemma~\ref{de0} and the limit~\eqref{thm:ld}, we have
  \begin{multline*}
   \lim_{V\to\infty}\sup_{m\ge y(qV)}\sup_{A\in\cal{B}(\cup_{n\ge 1}(\N^*)^n)}  \P\left({\rm{bulk}}(X_1,\dots,X_{N(qV)})\in A,\left|N(qV)-qV\right|\ge \delta_0qV\bigg|\sum_{i=1}^{N(qV)}X_i=m\right)\\
    \le \lim_{V\to\infty}\sup_{m\ge y(qV)} \frac{\P\left(\left|N(qV)-qV\right|\ge \delta_0qV,\sum_{i=1}^{N(qV)}X_i=m\right)}{qV\P(X=\lceil m-\rho_c V\rceil )}\frac{qV\P(X=\lceil m-\rho_c V\rceil)}{\P\left(\sum_{i=1}^{N(qV)}X_i=m\right)}=0.
  \end{multline*}
  In conclusion, we have
   \begin{multline*}
    \limsup_{V\to\infty}\sup_{m\ge y(qV)}\sup_{A\in\cal{B}(\cup_{n\ge 1}(\N^*)^n)}\bigg|\P\left({\rm{bulk}}(X_1,\dots,X_{N(qV)})\in A\bigg|\sum_{i=1}^{N(qV)}X_i=m\right)\\-\P\left((X_1,\dots,X_{N(qV)-1})\in A\right)\bigg|\le 2\eps.
   \end{multline*}
   We obtain the result by letting $\eps\to 0$.
\end{proof}

\section{Proof of the condensation phenomenon}\label{sec:cpp}
In this section, we prove our main results Theorem~\ref{main:cond} and Theorem~\ref{main:clt} with the help of the compound Poisson process introduced in Section~\ref{sec:cond}.
\begin{proof}[Proof of Theorem~\ref{main:cond}]
By the weak law of large numbers, we have, for any $\eps>0$,
  \[
  \P\left(\left|\frac{1}{V}\sum_{i=1}^{N(qV)-1}X_i-\rho_c\right|\ge \eps\right)\to 0.
  \]
We concludes the limit~\eqref{lln} by using  Lemma~\ref{lmcpp2} and Proposition~\ref{poisum}.
Similarly, the limit~\ref{mgclt}  follows immediately from Lemma~\ref{lmcpp2}, Proposition~\ref{poisum} and the standard central limit theorems of Poisson random variables.
\end{proof}

\begin{proof}[Proof of Theorem~\ref{main:clt}]
By using the standard results of the compound Poisson processes and also Lemma~\ref{lmcpp2}, Proposition~\ref{poisum}, we can obtain the limits~\eqref{lln} and~\eqref{gauss} easily.

To prove the limit~\eqref{stable}, let
    \[
H(r):=\P\left(\sum_{i=1}^{N(q)}X_i>r\right)r^\alpha,
\]
then we apply the Theorem 3.37 in Foss et al.~\cite{foss} and find that the function $H$ is slowly varying, that is, for all $t>0$,
\[
\lim_{r\to\infty}\frac{H(tr)}{H(r)}=\lim_{r\to\infty}\frac{\P\left(\sum_{i=1}^{N(q)}X_i>tr\right)}{\P(X>tr)}\frac{L(tr)}{L(r)}\frac{\P(X>r)}{\P\left(\sum_{i=1}^{N(q)}X_i>tr\right)}=1.
\]
By using the standard central limit theorems for the sum of \emph{i.i.d.} random copies of $\sum_{i=1}^{N(q)}X_i$ (See Theorem 3.8.2. and Example 3.9.5. in Durrett~\cite{durrett}), Lemma~\ref{lmcpp2} and Proposition~\ref{poisum}, we obtain the result.

To prove the limit~\eqref{ga}, according to the proof Theorem~1 in Section 35 of the book Gnedenko and Kolmogorov~\cite{gnedenko1968limit}, there exists a sequence $C_V\to\infty$ such that
\begin{equation}\label{tt1}
V\P(X>C_V)\to0,
\end{equation}
\begin{equation}\label{tt2}
B_V^2:=V\sum_{j<C_V}j^2Q_j\phi_c^j\gg C_V^2,
\end{equation}
and
\begin{equation}\label{tt3}
\frac{(\sum_{j<C_V}jQ_j\phi_c^j)^2}{\sum_{j<C_V}j^2Q_j\phi_c^j}\to 0.
\end{equation}
Therefore, from the proof of the Lindeberg-Feller theorem (See for example the Theorem~3.4.10 in Durrett), we have
\[
V\left|\E\exp\{it(X\ind{X\le C_V}-\E X\ind{X\le C_V})/B_V\}-1+\frac{t^2}{2qV}\right|\to 0.
\]
It implies that
\[
\frac{\sum_{i=1}^{N(qV)}X_i\ind{X_i\le C_V}-N(qV)\E X\ind{X\le C_V}}{B_V}\Rightarrow \cal{N}(0,1)
\]

From~\eqref{tt1}, we see that
\[
\P\left(\sum_{i=1}^{N(qV)}X_i\neq \sum_{i=1}^{N(qV)}X_i\ind{X_i\le C_V}\right)\le qV\P\left(X>C_V\right)\to 0.
\]
Combining with
\[
\frac{(N(qV)-qV)\E X\ind{X\le C_V}}{\sqrt{V\sum_{j<C_V}j^2Q_j\phi_c^j}}\to 0.
\]
We obtain~\eqref{ga}.
\end{proof}

\section*{Acknowledgments}
  The author is supported by the National Key R\&D Program of China under Grant 2022YFA 1006500.

\bibliographystyle{amsplain}
\bibliography{ref}

\providecommand{\bysame}{\leavevmode\hbox to3em{\hrulefill}\thinspace}
\providecommand{\MR}{\relax\ifhmode\unskip\space\fi MR }
\providecommand{\MRhref}[2]{%
  \href{http://www.ams.org/mathscinet-getitem?mr=#1}{#2}
}
\providecommand{\href}[2]{#2}
\begin{thebibliography}{10}

\bibitem{adams2007large}
Stefan Adams, \emph{Large deviations for empirical path measures in cycles of
  integer partitions}, arXiv preprint math/0702053 (2007).

\bibitem{adams2008large}
\bysame, \emph{Large deviations for empirical cycle counts of integer
  partitions and their relation to systems of bosons}, Analysis and stochastics
  of growth processes and interface models (2008), 148--172.

\bibitem{aldous1999deterministic}
David~J Aldous, \emph{Deterministic and stochastic models for coalescence
  (aggregation and coagulation): a review of the mean-field theory for
  probabilists}, Bernoulli (1999), 3--48.

\bibitem{andreis2021large}
Luisa Andreis, Wolfgang K{\"o}nig, and Robert~IA Patterson, \emph{A
  large-deviations principle for all the cluster sizes of a sparse
  erd{\H{o}}s--r{\'e}nyi graph}, Random Structures \& Algorithms \textbf{59}
  (2021), no.~4, 522--553.

\bibitem{armen}
In{\'e}s Armend{\'a}riz and Michail Loulakis, \emph{Thermodynamic limit for the
  invariant measures in supercritical zero range processes}, Probability theory
  and related fields \textbf{145} (2009), no.~1, 175--188.

\bibitem{armendariz}
In{\'e}s Armend{\'a}riz and Michail Loulakis, \emph{Conditional distribution of
  heavy tailed random variables on large deviations of their sum}, Stochastic
  Processes and their Applications \textbf{121} (2011), no.~5, 1138--1147.

\bibitem{ball}
{J. M.} Ball, J.~Carr, and O.~Penrose, \emph{The becker-d{\"o}ring cluster
  equations: Basic properties and asymptotic behaviour of solutions},
  Communications in Mathematical Physics \textbf{104} (1986), no.~4, 657--692.

\bibitem{MR1429811}
A.~Baltrunas, \emph{On a local limit theorem on one-sided large deviations for
  dominated-variation distributions}, Liet. Mat. Rink. \textbf{36} (1996),
  no.~1, 1--9. \MR{1429811}

\bibitem{betz}
Volker Betz, Daniel Ueltschi, and Yvan Velenik, \emph{Random permutations with
  cycle weights}, Ann. Appl. Probab. \textbf{21} (2011), no.~1, 312--331.
  \MR{2759204}

\bibitem{cline}
Daren B.~H. Cline, \emph{Intermediate regular and {$\Pi$} variation}, Proc.
  London Math. Soc. (3) \textbf{68} (1994), no.~3, 594--616. \MR{1262310}

\bibitem{dvz}
A.~Dembo, A.~Vershik, and O.~Zeitouni, \emph{Large deviations for integer
  partitions}, Markov Process. Related Fields \textbf{6} (2000), no.~2,
  147--179. \MR{1778750}

\bibitem{dembo}
Amir Dembo and Ofer Zeitouni, \emph{Large deviations techniques and
  applications}, Springer, 1998.

\bibitem{denisov}
D.~Denisov, A.~B. Dieker, and V.~Shneer, \emph{Large deviations for random
  walks under subexponentiality: the big-jump domain}, Ann. Probab. \textbf{36}
  (2008), no.~5, 1946--1991. \MR{2440928}

\bibitem{durrett}
Richard Durrett, \emph{Probability: theory and examples}, second ed., Duxbury
  Press, Belmont, CA, 1996. \MR{MR1609153 (98m:60001)}

\bibitem{Durrett99}
Richard Durrett, Boris~L. Granovsky, and Shay Gueron, \emph{The equilibrium
  behavior of reversible coagulation-fragmentation processes}, J. Theoret.
  Probab (1999), 447--474.

\bibitem{foss}
Sergey Foss, Dmitry Korshunov, Stan Zachary, et~al., \emph{An introduction to
  heavy-tailed and subexponential distributions}, vol.~6, Springer, 2011.

\bibitem{gnedenko1968limit}
B.~V. Gnedenko and A.~N. Kolmogorov, \emph{Limit distributions for sums of
  independent random variables}, Translated from the Russian, annotated, and
  revised by K. L. Chung. With appendices by J. L. Doob and P. L. Hsu. Revised
  edition, Addison-Wesley Publishing Co., Reading, Mass.-London-Don Mills.,
  Ont., 1968. \MR{MR0233400 (38 \#1722)}

\bibitem{gross}
Stefan Gro{\ss}kinsky, Gunter~M. Sch\"{u}tz, and Herbert Spohn,
  \emph{Condensation in the zero range process: stationary and dynamical
  properties}, J. Statist. Phys. \textbf{113} (2003), no.~3-4, 389--410.
  \MR{2013129}

\bibitem{hingant2017deterministic}
Erwan Hingant and Romain Yvinec, \emph{Deterministic and stochastic
  becker--d{\"o}ring equations: Past and recent mathematical developments},
  Stochastic processes, multiscale modeling, and numerical methods for
  computational cellular biology, Springer, 2017, pp.~175--204.

\bibitem{jansen2015large}
Sabine Jansen, Wolfgang K{\"o}nig, and Bernd Metzger, \emph{Large deviations
  for cluster size distributions in a continuous classical many-body system},
  The Annals of Applied Probability \textbf{25} (2015), no.~2, 930--973.

\bibitem{jeon1998existence}
Intae Jeon, \emph{Existence of gelling solutions for coagulation-fragmentation
  equations}, Communications in mathematical physics \textbf{194} (1998),
  no.~3, 541--567.

\bibitem{jeon2}
Intae Jeon, Peter March, and Boris Pittel, \emph{Size of the largest cluster
  under zero-range invariant measures}, Ann. Probab. \textbf{28} (2000), no.~3,
  1162--1194. \MR{1797308}

\bibitem{kelly}
F.P. Kelly, \emph{Reversibility and stochastic networks}, Probability and
  Statistics Series, J. Wiley, 1979.

\bibitem{kumi}
C.~Kl\"{u}ppelberg and T.~Mikosch, \emph{Large deviations of heavy-tailed
  random sums with applications in insurance and finance}, J. Appl. Probab.
  \textbf{34} (1997), no.~2, 293--308. \MR{1447336}

\bibitem{Mayer}
J.E. Mayer, \emph{Statistical mechanics}, Wiley, 1950.

\bibitem{MR0282396}
A.~V. Nagaev, \emph{Limit theorems that take into account large deviations when
  {C}ram\'{e}r's condition is violated}, Izv. Akad. Nauk UzSSR Ser. Fiz.-Mat.
  Nauk \textbf{13} (1969), no.~6, 17--22. \MR{282396}

\bibitem{sun2018functional}
Wen Sun, \emph{A functional central limit theorem for the becker--d{\"o}ring
  model}, Journal of Statistical Physics \textbf{171} (2018), no.~1, 145--165.

\bibitem{sunmdp}
\bysame, \emph{{A conditional compound Poisson process approach to the sparse
  Erd{\H{o}}s-R{\'e}nyi random graphs: moderate deviations}},  (2023),
  arXiv:2310.06348.

\bibitem{ldpsun}
\bysame, \emph{Pathwise large deviations for the pure jump k-nary interacting
  particle systems}, Ann. Appl. Probab. \textbf{34} (2024), no.~1A, 743--794.

\bibitem{MR0368115}
S.~G. Tka\v{c}uk, \emph{Local limit theorems, allowing for large deviations, in
  the case of stable limit laws}, Izv. Akad. Nauk UzSSR Ser. Fiz.-Mat. Nauk
  \textbf{17} (1973), no.~2, 30--33, 70. \MR{368115}

\bibitem{vershik}
A.~M. Vershik, \emph{Statistical mechanics of combinatorial partitions, and
  their limit configurations}, Funktsional. Anal. i Prilozhen. \textbf{30}
  (1996), no.~2, 19--39, 96. \MR{1402079}

\end{thebibliography}

\end{document}